\documentclass[12pt]{article}
\usepackage{amssymb,amsmath}
\newtheorem{thm}{Theorem}
\newtheorem{cor}[thm]{Corollary}
\newtheorem{lemma}[thm]{Lemma}
\newenvironment{defin}{\medskip\noindent{\sc
Definition}.}{\goodbreak\medskip}
\newenvironment{nota}{\medskip\noindent{\sc
Notations}.}{\goodbreak\medskip}
\newenvironment{remk}{\noindent{\sc
Remark}.}{\goodbreak\vskip10pt}
\newtheorem{prop}[thm]{Proposition}

\def\demo{\medskip\goodbreak\noindent
     \hbox{\sc Proof \kern .3em}\ignorespaces}%
  \def \qedbox{$\square$}%
  \def \qed{\hglue1mm\hfill{\ifmmode\qedbox
     \else\unskip\ \hglue0mm\hfill\qedbox\medskip
      \goodbreak\fi}}%
\def\enddemo{\qed\goodbreak\vskip10pt}%

\def\qed{\hglue1mm\hfill\raise -2pt\hbox{\vrule\vbox to 10pt{\hrule width
4pt
                  \vfill\hrule}\vrule}}

\newcommand{\T}{\mathbb {T}}
\newcommand{\A}{\mathbb {A}}

\newcommand{\R}{\mathbb {R}}

\newcommand{\Z}{\mathbb {Z}}
\newcommand{\N}{\mathbb {N}}

\newcommand{\Vc}{\mathcal {V}}

\newcommand{\Pc}{\mathcal {P}}

\newcommand{\Cc}{\mathcal {C}}
\newcommand{\Ic}{\mathcal {I}}
\newcommand{\Rc}{\mathcal {R}}
\newcommand{\Mc}{\mathcal {M}}

\newcommand{\Gc}{\mathcal {G}}
\newcommand{\Lc}{\mathcal {L}}

\newcommand{\Tc}{\mathcal {T}}

\begin{document}
\title{{
}
The link between the shape of the Aubry-Mather sets and their Lyapunov exponents
 }
\author{M.-C. ARNAUD
\thanks{ANR KAM faible}
\thanks{Universit\'e d'Avignon et des Pays de Vaucluse, Laboratoire d'Analyse non lin\' eaire et G\' eom\' etrie (EA 2151),  F-84 018Avignon,
France. e-mail: Marie-Claude.Arnaud@univ-avignon.fr}
}
\maketitle
\abstract{  We consider the irrational Aubry-Mather sets of an exact symplectic monotone $C^1$ twist map, introduce  for them a notion of ``$C^1$-regularity'' (related to the notion of Bouligand paratingent cone) and prove that~:
\begin{enumerate}
\item[$\bullet$] a Mather measure has zero Lyapunov exponents iff its support is almost everywhere $C^1$-regular;
\item[$\bullet$] a Mather measure has non zero Lyapunov exponents iff its support is almost everywhere $C^1$-irregular;
\item[$\bullet$] an Aubry-Mather set is uniformly hyperbolic iff it is everywhere non regular;
\item[$\bullet$] the Aubry-Mather sets which are close to the KAM invariant curves, even if they may be non $C^1$-regular, are not ``too irregular'' (i.e. have small paratingent cones).
\end{enumerate}
The main tools that  we use in the proofs are the so-called Green bundles.

}
\newpage
\tableofcontents
\newpage
 
\section{Introduction}

The  exact symplectic twist maps were studied for a long time because they represent (via a
symplectic change of coordinates) the dynamic of the generic symplectic diffeomorphisms of
surfaces near their elliptic periodic points (see \cite{Ch1}). One motivating  example of such  a map was introduced by Poincar\'e for the study of  the restricted 3-Body problem.\\

\noindent For these maps, the first invariant sets which were studied were the periodic orbits~: the ``last geometric Poincar\'e's theorem'' was proved by G.~D.~Birkhoff in 1913 in \cite{Bir2}. Later, in the 50's, the K.A.M.~theorems provide the existence of some invariant curves  for sufficiently regular symplectic diffeomorphisms of surfaces near their elliptic fixed points (see \cite{Ko}, \cite{Arno}, \cite{Mo} and \cite{Ru}). Then, in the 80's, the Aubry-Mather sets were discovered simultaneously and independently by Aubry \& Le Daeron (in \cite{ALD}) and Mather (in \cite{Mat1}). These sets are the union of  some quasi-periodic (in a weak sense) orbits, which are not necessarily on an invariant curve. We can define for each of these sets a {\em rotation number} and for every real number, there exists at least one Aubry-Mather set with this rotation number.\\ 

\noindent In 1988, Le Calvez proved in \cite{LC2} that for every generic exact symplectic   twist map $f$, there exists an open dense subset $U(f)$ of $\R$ such that every Aubry-Mather set for $f$ whose rotation number belongs to $U(f)$ is hyperbolic. Of course it doesn't imply that {\em all} the Aubry-Mather sets are hyperbolic (in particular    the K.A.M.~curves are not hyperbolic). \\
Some results are known concerning these hyperbolic Aubry-Mather sets~: it is proved in \cite{MMP}  that their projections have zero Lebesgue measure and in \cite{Mac1} that they have zero Hausdorff dimension.

\noindent The main question which will interest ourselves is then~: given some Aubry-Mather set of a symplectic twist map, is there a link between the geometric shape of these set and the fact that it is hyperbolic? Or ~: can we deduce the Lyapunov exponents of the measure supported on the Aubry-Mather set from the ``shape'' of this measure?\\
I didn't hear of such results for any dynamical systems and I think that the ones contained in this article are the first in this direction.
\\

\noindent Before explaining what kind of positive answers we can give to this question, let us introduce some notations and definitions. For classical results concerning exact symplectic twist map, the reader is referred to the books \cite{Go1} or \cite{Lc1}.\\

\begin{nota}
\noindent $\bullet$ $\T=\R/\Z$ is the circle.

\noindent $\bullet$ $\A=\T\times \R$ is the annulus and an element of $\A$ is denoted
by $(\theta, r)$.

\noindent $\bullet$ $\A$ is endowed with its usual symplectic form, $\omega=d\theta\wedge dr$ and its usual Riemannian metric.

\noindent  $\bullet$ $\pi~: \T \times \R \rightarrow\T$ is the first projection and $\tilde\pi~:
\R^2\rightarrow \R$ its lift. 

\noindent $\bullet$ $p~: \R^2\rightarrow \A$ is the usual covering map.

\end{nota}

\begin{defin} A $C^1$ diffeomorphism $f~: \A\rightarrow \A$ of the annulus which is isotopic
to identity  is a {\em positive twist map} if, for any given lift $\tilde f~: \R^2\rightarrow
\R^2$ and for every
$\tilde\theta\in\R$, the maps $r\mapsto \tilde\pi\circ \tilde f(\tilde\theta,r)$ and $r\mapsto \tilde\pi\circ \tilde
f^{-1}(\tilde\theta,r)$ are both diffeomorphisms, the first one increasing and the second one
decreasing. If $f$ is a positive twist map, $f^{-1}$ is a {\em negative} twist map. A {\em
twist map} may be positive or negative.\\ Moreover,
$f$ is {\em exact symplectic} if the 1-form
$f^*(rd\theta)-rd\theta$ is exact.

\end{defin}

\begin{nota} $\Mc_\omega^+$ is the set of  exact symplectic positive $C^1$ twist maps of $\A$,
$\Mc_\omega^-$ is the set of  exact symplectic negative $C^1$ twist maps of $\A$ and
$\Mc_\omega=\Mc_\omega^+\cup\Mc_\omega^-$ is the set of exact symplectic $C^1$ twist maps of
$\A$.

\end{nota}

\begin{defin} Let $M$ be a non-empty subset  of $\A$, let  $f~: \A\rightarrow \A$ be an exact symplectic twist map and let $\tilde f~:\R^2\rightarrow \R^2$ be one of its lifts. The set $M$ is {\em $f$-ordered} if~:
\begin{enumerate}
\item[$\bullet$] $M$ is compact;
\item[$\bullet$] $M$ is $f$-invariant;
\item[$\bullet$] $\forall z,z'\in p^{-1}(M), \tilde\pi(z)<\tilde\pi (z')\Leftrightarrow \tilde\pi(\tilde f (z))<\tilde\pi(\tilde f (z'))$
\end{enumerate}
(let us notice that this definition doesn't depend on the choice of the lift $\tilde f$ of $f$).

\end{defin}
A classical result asserts that every $f$-ordered set is a Lipschitz graph above a compact part of the circle. Moreover, if $K$ is a compact part or $\A$, there exists a constant $k>0$ depending only on $K$ and $f$ such that the Lipschitz constant of every $f$-ordered set meeting $K$ is less than $k$.

\begin{defin}
An {\em Aubry-Mather set} for an exact symplectic twist map $f$ is a minimal (for ``$\subset$'') $f$-ordered set. 
\end{defin}
Then it is well-known that if $M$ is an Aubry-Mather set of a $f\in \Mc_\omega$, there exists a bi-Lipschitz orientation preserving homeomorphism $h~:\T\rightarrow \T$ of the circle such that~: $\forall (\theta, r)\in M, \pi\circ f(\theta , r)=h(\theta)$~: the dynamic of $f$ on $M$ is conjugate via the first projection to the one of a bi-Lipschitz homeomorphism of the circle on a minimal invariant compact set. If we write the previous equality for a lift $\tilde f$ of $f$, we can associate to every Aubry-Mather set $M$ of $f$ a rotation number (which is the rotation number of any $\tilde h$ such that~: $\forall (\tilde \theta, r)\in \widetilde{M}=p^{-1}(M), \tilde h(\tilde \theta)=\tilde\pi\circ\tilde f(\tilde\theta, r)$) denoted by $\rho (M, \tilde f)$. Then for every $\rho\in\R$, there exists at least one Aubry-Mather set $M$ for $f$ such that $\rho(M, \tilde f)=\rho$. With our definition of Aubry-Mather set (minimal), if $\rho(M,\tilde f)$ is rational, then $M$ is a periodic orbit; in the other case, we will say that the Aubry-Mather set is {\em irrational} and two cases may happen~:
\begin{enumerate}
\item[$\bullet$] either $M$ is a curve (and $h$ is $C^0$-conjugate to a rotation);
\item[$\bullet$] or $M$ is a Cantor (and $h$ is a Denjoy counter-example).
\end{enumerate}
Moreover, every Aubry-Mather set carries a unique $f$-invariant Borel probability measure, denoted by $\mu(M,f)$. This measure is always ergodic (even uniquely ergodic on its support) and its support is $M$. Such a measure $\mu$ (associated to an Aubry-Mather set $M$ for $f$) will be called a {\em Mather measure}. 

Let us now explain what we mean by ``shape of a set'' or of a measure. This notion is related to a notion of regularity~:

\begin{defin}
Let $M\subset \A$ be a subset of $\A$ and $x\in M$ a   point of $M$. The   {\em paratingent cone} to $M$ at $x$ is the cone of $T_x\A$ denoted by $P_M(x)$ whose elements  are the limits~:
$$v=\lim_{ n\rightarrow \infty} \frac{x_n-y_n}{ t_n} $$
where $(x_n)$ and $(y_n)$ are sequences of elements of $M$ converging to $x$, $(t_n)$ is a sequence of elements of $\R_+^*$ converging to $0$, and $x_n-y_n\in\R$, refers to the unique lift of this element of $\A$ which belongs to $[-\frac{1}{2}, \frac{1}{2}[^2$.\\
\noindent We will say that $M$ is {\em $C^1$-regular} at $x$ if there exists a line $D$ of $T_x\A$ such that $P_M(x)\subset D$.
\end{defin}
This notion of (Bouligand's) paratingent cone comes from non-smooth analysis (see for example \cite{A-F}). Of course, at an isolated point, the notion of regularity doesn't mean anything, and we will use it only for Aubry-Mather sets having no isolated point, i.e. irrational Aubry-Mather sets.

\begin{thm}\label{th1}
Let $f\in\Mc_\omega$ be an exact  symplectic twist map and let $\mu$ be an irrational Mather measure of $f$.  The two following assertions are equivalent~:
\begin{enumerate}
\item[$\bullet$] for $\mu$-almost every $x$, ${\rm supp}\mu$ is $C^1$-regular at $x$;
\item[$\bullet$] the Lyapunov exponents of $\mu$ (for $f$) are   zero.
\end{enumerate}
\end{thm}
An alternative statement of this result is~: 

\begin{prop}\label{prop2}
Let $f\in\Mc_\omega$ be an exact  symplectic twist map and let $\mu$ be an irrational Mather measure of $f$.  The two following assertions are equivalent~:
\begin{enumerate}
\item[$\bullet$] for $\mu$-almost every $x$, ${\rm supp}\mu$ is not $C^1$-regular at $x$;
\item[$\bullet$] the Lyapunov exponents of $\mu$ (for $f$) are  non-zero.
\end{enumerate}
\end{prop}
Hence we don't obtain exactly the kind of result we wanted~: knowing the measure $\mu$ (and not the diffeomorphism $f$!), we can say if the Lyapunov exponents are zero or not, but the a priori knowledge of the Aubry-Mather set is not sufficient to deduce if the Lyapunov exponents are zero or no. To precise this fact, it would be interesting to answer to the following questions~:  
\medskip

\noindent{\bf Questions~:} \\
$\bullet$ Let us assume that $M$ is an irrational Aubry-Mather set of an exact symplectic $C^1
$ twist map $f$. Does there exist another exact symplectic  $C^1$ twist map $g$ such that $M$ is an irrational Aubry-Mather set for $g$ and such that $\mu (M, f)$ and $\mu (M, g)$ are not equivalent (i.e. not mutually absolutely continuous)?\\
$\bullet$ Does there exist an  irrational Aubry-Mather set $M\subset \A$ of an exact symplectic $C^1
$ twist map $f$, such that for every exact symplectic $C^1
$ twist map $g$ for which $M$ is an irrational Aubry-Mather set, the measures $\mu (M,f)$ and $\mu(M,g)$ are equivalent?

\medskip

However, in the extreme cases, we obtain   a result concerning the shape of the Aubry-Mather sets~:

 \begin{cor}\label{Cor3}
Let $f\in\Mc_\omega$ be an exact  symplectic twist map and let $M$ be an irrational Aubry-Mather set of $f$. If for all $x\in M$, $M$ is $C^1$-regular at $x$, then the Lyapunov exponents of $\mu(M, f)$ (for $f$) are zero.
\end{cor}

It is not hard to see that an Aubry-Mather set is everywhere $C^1$-regular if and only if there exists a $C^1$ map $\gamma~: \T\rightarrow \R$ whose graph contains $M$. In \cite{He2}, M.~Herman gives some examples of irrational Aubry-Mather sets which are invariant by a twist map, contained in a $C^1$-graph but not contained in an {\em invariant} continuous curve. I don't know any example of an irrational Aubry-Mather set with zero Lyapunov exponents which is not contained in a $C^1$ curve.

\medskip
\noindent{\bf Problem~:} is it possible to build an irrational Aubry-Mather set with zero Lyapunov exponents which is not contained in a $C^1$ graph?

\begin{prop}\label{prop4} Let $f\in\Mc_\omega$ be an exact  symplectic twist map and let $M$ be an irrational Aubry-Mather set of $f$.  The two following assertions are equivalent~:
\begin{enumerate}
\item[$\bullet$] for all $x\in M$, $M$ is not $C^1$-regular at $x$;
\item[$\bullet$] the set $M$ is uniformly hyperbolic (for $f$).
\end{enumerate}

\end{prop}\label{prop5}
In the non uniformly hyperbolic case, we can be more specific~: 

\begin{prop} \label{prop5}Let $f\in\Mc_\omega$ be an exact  symplectic twist map and let $\mu$ be an irrational  Mather measure of $f$ which is non uniformly hyperbolic, i.e. the Lyapunov exponents are non zero but the corresponding Aubry-Mather set $M={\rm supp}\mu$ is not (uniformly) hyperbolic. 
Then there exists a dense $G_\delta$ subset $\Gc$ of $M$ such that $M$ is $C^1$-regular at every point of $\Gc$.

\end{prop}
I must say that I don't know any example of an irrational Aubry-Mather set which is non uniformly hyperbolic.

\medskip

Let us now consider   
what happens near a K.A.M. invariant 
curve $\Cc$ for a generic $f\in \Mc_\omega$~: 
if $\alpha$ is the rotation number of this K.A.M. 
curve, for every neighbourhood $\Vc$ of $\Cc$ for 
the Hausdorff topology, there exists $\varepsilon >0$ 
such that every Aubry-Mather set whose rotation 
number is in $]\alpha-\varepsilon, \alpha+\varepsilon[$
 belongs to $\Vc$ (indeed, a limit of $f$-ordered set
 is $f$-ordered and the rotation number is continuous
 on the set of $f$-ordered sets; moreover, a classical
 result asserts that if there is a KAM curve, it is 
the unique $f$-ordered set having this rotation number).
 Hence, using Le Calvez' result mentioned before, we
 find in every neighbourhood $\Vc$ of $\Cc$ some 
irrational  uniformly hyperbolic Aubry-Mather sets, 
and hence some $C^1$-irregular Cantor sets (see the beginning
of the proof of proposition \ref{prop30} to see why
it cannot be a curve). But even if these Cantor sets
are
$C^1$-irregular, the closest they are to $\Cc$, the
less irregular they are  in the following sense~:

\begin{thm}\label{th6}
Let $f\in\Mc_\omega$ be an exact symplectic twist map and $\Cc$ be a $C^1$ invariant curve which is a graph such that $f_{|\Cc}$ is $C^1$ conjugate to a rotation. Let $W$ be a neighbourhood  of $T^1\Cc$, the unitary tangent bundle to $\Cc$  in  $T^1\A$, the unitary tangent bundle to $\A$. Then there exists a neighbourhood $V$ of $\Cc$ in $\A$ such that for every Aubry-Mather set $M$ for $f$ contained in $V$~:
$$\forall x\in M, P^1_M(x)\subset W$$
where $P^1_M(x)$ refers to the unitary paratingent cone.
\end{thm}
It implies that in this case, even if the paratingent cone at $x$ to $M$ is not a line, it is a thin cone close to a line.\\

To prove the results contained in this article, we will use a very usefull mathematical object~: the Green bundles. 
They were introduced by L.~W.~Green in \cite{Gr} for Riemannian
geodesic flows; then P.~Foulon extended this construction to Finsler metrics in \cite{Fo1}
and  G.~Contreras and  R.~Iturriaga extended it in
\cite{C-I} to optical Hamiltonian flows; in \cite{Bi-Ma}, M.~Bialy and R.~S.~Mackay give an
analogous construction for the dynamics of sequence of symplectic twist maps of $T^*\T ^d$
without conjugate point. Let us cite also a very short survey \cite{It1} of R.~Iturriaga  on
the various uses of these bundles (problems of rigidity, measure of hyperbolicity\dots). \\

In \cite{Arna1} and \cite{Arna2}, I constructed these bundles along invariant graphs and proved, under various dynamical assumptions, that they may be used to prove some results of $C^1$- regularity. In particular, the strongest result contained in \cite{Arna1} for twist maps is  that the ``Birkhoff invariant curves'' are more regular than Lipschitz
(more precisely $C^1$ regular on a dense $G_\delta$
subset) or, equivalently  that the
$C^1$ solutions of the Hamilton-Jacobi equation are
Lebesgue almost everywhere $C^2$.\\

In the second section of this new article, I enlarge the construction of the Green bundles to the Aubry-Mather sets, give some of their properties (semi-continuity\dots), introduced a notion of $C^1$-regularity (which is quite different from the one contained in \cite{Arna1}) and explain how the coincidence of the two Green bundles implies some regularity of the Aubry-Mather sets.\\
In the third section, I   explain how the (almost everywhere) transversality of the Green bundles implies some (non uniform) hyperbolicity. This result concerning Lyapunov exponents is completely new. In the case of uniform hyperbolicity, it is a consequence of a result of Contreras and Ituriaga, but we prove even in this case a more precise result (we don't assume that the dynamic is non wandering). We recall some well known results too.\\
In the fourth section, we prove that hyperbolic Aubry-Mather sets are $C^1$-irregular. These results too are completely new, and we deal with the uniformly and non uniformly hyperbolic cases.\\
 Finally, in the last section, we prove the results contained in the introduction.

\medskip
\noindent{\small\sc Acknowledgments}. {\small I am  grateful  to R.~Perez-Marco who firstly suggested me that the result for Aubry-Mather sets could be ``hyperbolicity versus regularity'', to J.-C.~Yoccoz whose questions   led me to the appropriate  definition of regularity,  to L.~Rifford   who pointed to
me the notion of Bouligand's paratingent cone and to S.~Crovisier who suggested me to send one Green bundle on the ``horizontal'' for the proof of the ``dynamical criterion'' , which gives a significant improvement of the proof.}

\section{Construction of the Green bundles along an irrational Aubry-Mather set, link with the $C^1$-regularity}

\begin{nota} $\pi~: \T \times \R \rightarrow\T$ is the projection.

\noindent If $x\in \A $, $V(x)=\ker D\pi (x)\subset T_x\A $ is the {\em vertical} at $x$.

\noindent If $x\in \A$ and $k\in \Z^*$, $G _ k(x)=Df^k(f^{-k}(x))V(f^{-k}(x))$is a 1-dimensional linear subspace (or line) of
$T_x\A
$.
\end{nota}

\begin{defin}
 If we identify $T_x\A$ with $\R^2$ by using  the standard coordinates 
$(\theta, r)\in\R^2$, we
may deal with the {\em slope} $s(L)$ of any line $L$ of $T_x\A$ which is transverse to the
vertical
$V(x)$~: it means that $L=\{ (t, s(L)t); t\in\R\}$.

  If $x\in \A$ and if $L_1$, $L_2$ are two lines of $T_x\A$ which are transverse to
the vertical
$V(x)$, $L_2$ is {\em above} (resp. {\em strictly above}) $L_1$ if
$s(L_2)\geq s(L_1)$ (resp. $s(L_2)>s(L_1)$). In this case, we write~: $L_1\preceq L_2$ (resp. 
 $L_1\prec L_2$). In a similar way, if $\Lc_1$ and $\Lc_2$ are two sets of lines of  $T_x\A$ which are transverse to
the vertical
$V(x)$, $\Lc_2$ is {\em above} (resp. {\em strictly above}) $\Lc_1$ if
$s(L_2)\geq s(L_1)$ (resp. $s(L_2)>s(L_1)$) for all $L_1\in \Lc_1, L_2\in \Lc_2$. In this case, we write~: $\Lc_1\preceq \Lc_2$ (resp. 
 $\Lc_1\prec \Lc_2$).

 A sequence $(L_n)_{n\in\N}$ of  lines of $T_x\A$ is {\em non decreasing}
(resp. {\em increasing}) if for every $n\in\N$, $L_n$ is transverse to the vertical and 
$L_{n+1}$ is above (resp. strictly above) $L_n$. We define the {\em non increasing} and
{\em decreasing} sequences of  lines of
$T_xM$ in a similar way.

\end{defin}
\begin{remk}
A decreasing sequence of lines corresponds to a decreasing sequence of slopes.
\end{remk}

\begin{defin}
If $K$ is a subset of $\A$ or of its universal covering $\R\times\R$, if $F$ is a
1-dimensional sub-bundle of $T_K\A$ (resp. $T_K\R^2$) transverse to the  vertical, we say that
$F$ is upper (resp. lower) semi-continuous if the map which maps $x\in K$ onto the slope
$s(F(x))$ of
$F(x)$ is upper (resp. lower) semi-continuous.
\end{defin}

\begin{prop}\label{prop7}
Let $f~: \T\times \R\rightarrow \T\times\R$ be an exact symplectic positive $C^1$ twist map
and let $M$ be a $f$-ordered  set. Then , for every $x\in M$ which is not an isolated point of $M$,  we have~: 
$$\forall n\in\N^*, G_{-n}(x)\prec G_{-(n+1)}(x)\prec P_M(x)\prec G _{n+1}(x)\prec G _n(x).$$

\end{prop}
(in this statement we identify the cone $P_M(x)$ with the set of the lines which are contained in this cone)\\
As an irrational Aubry-Mather set has no isolated point, we deduce~:
\begin{cor}
Let $f~: \T\times \R\rightarrow \T\times\R$ be an exact symplectic positive $C^1$ twist map
and let $M$ be an irrational Aubry-Mather set of $f$. Then , for every $x\in M$,  we have~: 
$$\forall n\in\N^*, G_{-n}(x)\prec G_{-(n+1)}(x)\prec P_M(x)\prec G _{n+1}(x)\prec G _n(x).$$

\end{cor}

\noindent{\bf Proof of proposition \ref{prop7}~:} As $M$ is a $f$-ordered set, it is the graph of a Lipschitz map $\gamma$ above a non empty and compact part $K$ of $\T$.  Let now $x=(t,\gamma(t))$ be a point of $M$.
We will use the left and right paratingent cones to $M$ at $x$, defined by~: 
\begin{enumerate}
\item[$\bullet$] the right paratingent cone of $M$ at $x$, denoted by $P_M^r(x)$, is the set whose elements  are the limits~:
 $\displaystyle{v=\lim_{ n\rightarrow \infty} \frac{(u_n, \gamma(u_n))-(s_n,\gamma(s_n))}{ t_n} }$ 
where $(u_n)$ and $(s_n)$ are sequences of elements of $K$ converging to $t$ from above (i.e. $u_n, s_n\in [t, +\infty[$) and $(t_n)$ is a sequence of elements of $\R_+^*$ converging to $0$;
 \item[$\bullet$] similarly, the left paratingent cone of $M$ at $x$, denoted by $P_M^l(x)$, is the set whose elements  are the limits~:
 $\displaystyle{v=\lim_{ n\rightarrow \infty} \frac{(u_n, \gamma(u_n))-(s_n,\gamma(s_n))}{ t_n}} $ 
where $(u_n)$ and $(s_n)$ are sequences of elements of $K$ converging to $t$ from below and $(t_n)$ is a sequence of elements of $\R_+^*$ converging to $0$.
\end{enumerate}

It is not hard to verify that every element of $P_M(x)$ is in the convex hull of $P_M^l(x)\cup P_M^r(x)$ (we identify the lines of $T_x\A$ transverse to the vertical with their slopes in order to deal with their convex hull). Hence, we only have to prove the inequalities of  proposition \ref{prop7} for $P_M^r(x)$ and $P_M^l (x)$ (and even for those of these two cones which are not trivial) to deduce the inequalities of this proposition. 
Because the four proofs are similar, we will assume  for example that $P_M^r(x)\not=\{ 0\}$ and we will prove that~: $\forall n\in\N^*, P_M^r(x)\prec G_{n+1}(x)\prec G_n(x)$.

In fact we shall need to deal with half lines instead of lines. Hence we define $\Pc_M^r(x)$ as being the set of the half lines of $T_x\A$ which are contained in $P_M^r(x)$ such that their points have positive abscissa. Equivalently,  $\Pc_M^r(x)$ is the set of the limits~:
 $\displaystyle{v=\lim_{ n\rightarrow \infty} \frac{(u_n, \gamma(u_n))-(s_n,\gamma(s_n))}{ t_n} }$ 
where $(u_n)$ and $(s_n)$ are sequences of elements of $K$ converging to $t$ such that~: $\forall n, t\leq s_n<u_n$ and $(t_n)$ is a sequence of elements of $\R_+^*$ converging to $0$.
As $M$ is $f$-ordered, we have~: $\forall y\in M, Df(\Pc^r_M(y))=\Pc^r_M(y)$ (in particular the image through $Df$ of the right paratingent cone at $y$ is the right paratingent cone at the image $f(y)$). Hence~: $\forall k\in\Z , \Pc_M^r(f^kx)=Df^k(\Pc_M^r(x))$.  Let now $V_+(x)=\{ (0, R), R>0\}\subset T_x\A$ be the {\em upper vertical} at $x$ and let us denote by $g_k(x)$ the half line~: $g_k(x)=Df^k(f^{-k}(x))V_+(f^{-k}x)$.

Let us   look at the action of $Df$ on the half lines of the tangent linear spaces $T_{f^k(x)}\A$. As $f$ is a positive twist map, we have (identifying as before $T_{f(x)}\A$ with $\R^2$)~: $Df(x)(0,1)=(a,b)$ with $a>0$. If now $\R_+(\alpha,\beta)\in\Pc_M^r(x)$, we know that $\alpha>0$. Hence  the base $((\alpha, \beta), (0,1))$ is a direct base (for $\omega$) of $T_x\A$; as $Df(x)$ is symplectic, the image base $((\alpha', \beta'), (a,b))$ is direct too. It means exactly that the line $\R (a,b)=G_1(f(x))$ is strictly above the line $\R(\alpha ', \beta')$ of $P_M^r(f(x))$. Repeating this argument for every half line of $\Pc_M^r(x)$ and every point of the orbit of $x$, we obtain 
that~: $\forall k\in \Z, P_M^r(f^k(x))\prec G_1(f^k(x))$.

Let us consider the action of $Df$ on the circles bundle of the half lines along the orbit of $x$~: as $f$ is orientation preserving, this action preserves the orientation of the circles. Moreover, if these circles are oriented in the direct sense, then any half line of $\Pc_M^r(f^k(x))$, $g_1(f^k(x))$ and   $V_+(f^kx)$ are in the direct sense (let us recall that on the oriented circle, we can speak of the orientation of three points   but not of a pair). Hence their image under $Df$, $Df^2$, \dots are in the same order, i.e~: any half line of  $\Pc_M^r(f^k(x))$, 
$g_{n+1}(f^k(x))$ and  $g_n(f^k(x))$   are in the direct sense, and then~: $P_M^r(f^k(x))\prec G_{n+1}(f^k(x))\prec G_n(f^k(x))\prec\dots \prec G_1(f^k(x))$.
\enddemo

\begin{remk} Let us notice that  in the proof of proposition \ref{prop7}, we have seen that~:
$$\forall x\in M, \forall n\geq 1, D\pi\circ Df^n(x)(0,1)>0.$$
In a similar way, we have~:  $\forall x\in M, \forall n\geq 1, D\pi\circ Df^{-n}(x)(0,1)<0$.

\end{remk}

Hence $(G_n(x))$ is a strictly decreasing  sequence of lines of
$T_x\A$ which is   bounded  below. Then it tends to a limit $G_+(x)$. In a similar way,
the sequence $(G_{-n}(x))$ tends to a limit, $G_-(x)$.

\begin{defin} If $x\in \A$ belongs to an irrational Aubry-Mather set  $M$ of $f\in \Mc^+_\omega$,
the bundles
$G_-(x)$ and
$G_+(x)$ are called the {\em Green bundles} at $x$ associated to
$f$.

\end{defin}

\smallskip
\noindent{\bf Example}  Let us assume that $x\in M$ is a periodic hyperbolic periodic
point of $f$; then $G_+(x)=E^u(x)$ is the tangent space to the unstable manifold of $x$ and
$G_-(x)=E^s(x)$ is the tangent space to the stable manifold.

\smallskip
In fact, in order to build the Green bundles for $f$ at a point $x\in\A$, we don't need that $x$  belongs to 
a  $f$-ordered set. Let us introduce the exact
set which will be useful for us (the one along which
we can define the Green bundles)~:

\begin{defin}
Let $f\in\Mc_\omega^+(f)$ be a positive exact symplectic twist map. Then the {\em Green set} of $f$, denoted by $\Tc(f)$, is the sets of points $x\in \A$ such that~:
\begin{enumerate}
\item[$\bullet$] for all $n\geq 1$ and all $k\in\Z$, $D\pi\circ Df^n(f^kx)(0,1)>0$ and $D\pi\circ Df^{-n}(f^kx)(0,1)<0$;
\item[$\bullet$] or all $n\geq 1$ and all $k\in\Z$, \\
$G_{-n}(f^kx)=Df^{-n}(f^{n+k}x)V(f^{n+k}x)\prec Df^{-(n+1)}(f^{n+1+k}x)V(f^{n+1+k}x)  =$\\
$  G_{-(n+1)}(f^kx) \prec  G_{n+1}(f^kx) = Df^{n+1}(f^{-(n+1)+k}x)V(f^{-(n+1)+k}x) $\\
$\prec  Df^n(f^{-n+k}x)V(^{-n+k}x) = G_n(f^kx)  $.
\end{enumerate}
\end{defin}
Let us notice that the first point is not useful to define the Green bundles, but will be used in the next section to prove the so-called ``dynamical criterion''. Then we have~:

\begin{prop}\label{prop9}
Let $f\in \Mc_\omega^+$ be an exact symplectic $C^1$ positive twist map. Then $\Tc (f)$ is a non-empty, closed subset of $\A$ which contains every irrational Aubry-Mather set of $f$ and is invariant by $f$. At every $x\in \Tc(f)$, we can define $G_-(x)$ and $G_+(x)$.
\end{prop}
\begin{remk}
Let us notice that every essential invariant curve by $f\in \Mc^+_\omega$ is  a subset of $\Ic (f)$ (see \cite{Arna1}).
\end{remk}

\noindent{\bf Proof of proposition \ref{prop9}~:} The only things that we have to prove is that $\Tc (f)$ is  closed. \\
 Because $f$ is a positive twist map, we have for every $x\in \A$~: $D\pi\circ Df (x)(0,1)>0$ and $D\pi\circ Df^{-1}(x)(0,1)<0$. Hence for every $x\in \A$, $V(x)$ and $G_1(x)$ are transverse, and $V(x)$ and $G_{-1}(x)$ are transverse too. We deduce that for every $x\in \A$ and every $n\in\N^*$, $G_n(x)=Df^{-(n+1)}G_1(f^{n+1}x)$ and $G_{n+1}(x)=Df^{-(n+1)}V(f^{n+1}x)$ are transverse, and $G_{-(n+1)}(x)$ and $G_{-n}(x)$ are transverse.\\
Let us now consider $\Cc(f)$ the set of $x\in \A$ such that~:
\begin{enumerate}
\item[$\bullet$] for all $n\geq 1$, $D\pi\circ Df^n(x)(0,1)\geq 0$ and $D\pi\circ Df^{-n}(x)(0,1)\leq 0$;
\item[$\bullet$] for all $n\in\N^*$, $G_{-1}\preceq \dots \preceq G_{-n}(x) \preceq  G_{-(n+1)}(x)\preceq G_{n+1}(x)  \preceq  G_n(x) \preceq \dots \preceq G_1(x) $.
 \end{enumerate}
 Then $\Cc(f)$ is closed. If we prove that $\Cc (f)=\Tc (f)$, we have finished the proof.\\
 We have~: $\Tc (f)\subset \Cc(f)$. Moreover, if $x\in \Cc (f)$, we know that for all $n\in\N^*$, $G_{ n+1 }(x)\preceq G_{n }(x) $; as $G_{n}(x)$ and $G_{n+1}(x) $ are transverse, we deduce that $G_{n+1}(x)\prec G_n(x)$. In a similar way, we obtain that $G_{-n}(x) \prec  G_{-(n+1)}(x)$. From $G_{-n}(x) \prec  G_{-(n+1)}(x)\preceq G_{n+1}(x)\prec G_n(x)$, we deduce~: $G_{-n}(x)\prec G_n(x)$. Thus if $x\in\Cc(f)$, $x$ satisfies the second point of the definition of $\Tc (f)$. Hence every $G_k(x)$ for $k\in \Z^*$ is transverse to the vertical and~: $\forall k\in\Z^*, \forall x\in\Cc(f), D\pi\circ Df^k(x)(1,0)\not=0$. Therefore $x\in\Cc (f)$ satisfies the first point of the definition of $\Tc (f)$ too. Finally~: $\Cc(f)\subset\Tc (f)$ and then $\Cc(f)=\Tc (f)$.
\enddemo

Having built the Green bundles on $\Tc (f)$, we can give some of their properties, similar to the ones given in \cite{Arna1}, which in particular give a link between these Green bundles and the notion of $C^1$-regularity.

\begin{prop}\label{prop10} Let $f$ be  an exact symplectic positive $C^1$ twist map $f~: \A\rightarrow \A$. Then
the Green bundles, defined at every point of $\Tc (f)$, are invariant by $Df$.\\
The map $(x\in\Tc(f)\rightarrow G_+(x))$ is upper semi-continuous and the
map $x\rightarrow G_-(x)$ is lower semi-continuous and we have~: $\forall x\in \Ic (f), G_-(f)\preceq G_+(f)$. Therefore, the set~:
$$\Gc(f)=\{x\in \Tc (f); G_-(x)=G_+(x)\}$$ is a $G_\delta$ subset of $\Tc (f)$.\\
Moreover, for every irrational Aubry-Mather set $M$ of $f$ and every $x\in M$,  we have~: 
$G_-(x)\preceq P_M(x)\preceq G_+(x)$ and for
every $x_0\in  \Gc (f)\cap M$,  $M$ is $C^1$ regular  at $x_0$ and
$P_M(x_0)=G_+(x_0)=G_-(x_0)$. Moreover,
$G_-$ and $G_+$ are continuous at such a $x_0$.

\end{prop}

This proposition is a corollary of proposition \ref{prop7} and of usual properties of
real functions (the fact that the (simple) limit of a decreasing sequence of continuous
functions is upper semi-continuous).

\begin{cor}\label{C7} Let $M$ be an irrational Aubry-Mather set 
of an exact symplectic positive $C^1$ twist map $f~: \A\rightarrow \A$. 
  We
assume that~:
$$\forall x\in M, G_-(x)=G_+(x).$$
Then $M$ is $C^1$ regular at every $x\in M$ and there exists a $C^1$ map $\gamma~: \T\rightarrow \R$ whose graph contains $M$.\\
 Moreover, in this case, at every $x=(t,\gamma (t))\in M$, the sequences
$( G_n (x))_{n\in\N}$ and $( G_{-n}(x))_{n\in\N}$ converge uniformly to
$\R (1, \gamma'(t))$.
\end{cor}
Everything in this corollary is a consequence of proposition \ref{prop10}; the fact that the
convergence is uniform comes from Dini's theorem~: if an increasing or decreasing sequence of
real valued continuous functions defined on a compact set converges simply to a continuous
function, then the convergence is uniform.\\
This corollary gives us some criterion using the Green bundles to prove that an Aubry-Mather set is $C^1$-regular. But of course we never said that the transversality of the Green bundles implies the non regularity of the corresponding Aubry-Mather set. This will be explained later.
\section{Green bundles and Lyapunov exponents}

\subsection{A dynamical criterion}

We begin   by giving a criterion to determine if a given vector is in one of the two Green
bundles.

\begin{prop}\label{prop12}
Let $f$ be an exact symplectic positive $C^1$ twist map and let $x\in \Tc (f)$ be a point of the Green set  whose orbit $\{ f^k(x), k\in \Z\}$ is relatively compact.  Then~:
$$\lim_{n\rightarrow +\infty} D\pi\circ Df^n(x)(1,0)=+\infty\quad{\rm et}\quad \lim_{n\rightarrow +\infty} D\pi\circ Df^{-n}(x)(1,0)=-\infty .$$
\end{prop}

\begin{cor}\label{cor13} {\bf (dynamical criterion)}
Let $f$ be an exact symplectic positive $C^1$ twist map and let $x\in \Tc (f)$ be a point of the Green set  whose orbit $\{ f^k(x), k\in \Z\}$ is relatively compact.  Let $v\in T_x\A$. then~:
\begin{enumerate}
\item[$\bullet$] if $v\notin G_-(x)$ then~: $\displaystyle{\lim_{n\rightarrow +\infty} | D\pi\circ Df^n(x)v|=+\infty}$;
\item[$\bullet$] if $v\notin G_+(x)$ then~: $\displaystyle{\lim_{n\rightarrow +\infty} | D\pi\circ Df^{-n}(x)v|=+\infty}$.
\end{enumerate} 
\end{cor}

\noindent{\bf Proof of proposition \ref{prop12} and corollary \ref{cor13}~:} We will only prove the part of proposition and corollary corresponding to what happens in $+\infty$.
  We use the standard symplectic coordinates
$(\theta , r)$ of
$\A$ and we define for every $k\in \Z$~: $x_k=f^k(x)$. 
\noindent In these coordinates, for $j\in\Z^*$,  the line $G_j(x_k)$ is the graph of $(t\rightarrow
s_j(x_k)t)$ ($s_j(x_k)$ is the slope of $G_j(x_k)$). \\
The matrix $M_n(x_k)$ of $Df^n(x_k)$  (for $n\geq 1$)
is a symplectic matrix~:

$$M_n(x_k)=\begin{pmatrix}a_n(x_k)&b_n(x_k)\\
c_n(x_k)&d_n(x_k)\\
\end{pmatrix}$$
with $\det M_n(x_k)=1$. We know that the
coordinate $D(\pi\circ f^n)(x_k)(0, 1)=b_n(x_k)$ is strictly positive. Using the definition of
$G_n (x_{k+n})$, we obtain~: $d_n(x_k)=s_n (x_{k+n})b_n(x_k)$.\\
The matrix $M_n(x_k)$ being symplectic, we have~:
$$M_n(x_k)^{-1}=\begin{pmatrix}d_n(x_k)&-b_n(x_k)\\
-c_n(x_k)&a_n(x_k)\\
\end{pmatrix}$$
we deduce from the definition of $G_{-n}(x_k)$ that~: $a_n(x_k)=-b_n(x_k)s_{-n}(x_k)$. Finally,
if we use the fact that $\det M_n(x_k)=1$, we obtain~:
$$M_n(x_k)=\begin{pmatrix} -b_n(x_k)s_{-n}(x_k) & b_n(x_k)\\
 -b_n(x_k)^{-1}- b_n(x_k)s_{-n}(x_k)s_n (x_{k+n})    & s_n (x_{k+n})b_n(x_k)\\
\end{pmatrix}$$
\begin{lemma} \label{lem14}Let $K$ be a compact       subset of $\Tc(f)$. There exists a constant $A>0$ such that~:
$$\forall x\in K, \forall n\in\N^*, \max\{ |s_n (x)|, |s_{-n}(x)|\}\leq A.$$

\end{lemma}
\noindent{\bf Proof of lemma \ref{lem14}~:}
We deduce from the definition of $\Ic (f)$  that~: $\forall x\in \Tc (f) , \forall n\in\N^*,
s_{-1}(x)\leq s_{-n}(x) < s_n (x)\leq s_1 (x)$. Therefore, we only have to prove the
inequalities of the lemma for $n=1$. \\
The real number $s_{-1}(x)$, which is the slope of $Df^{-1}(f(x))V(f(x))$, depends
continuously on $x$, and is defined for every $x$ belonging to the compact subset $K$. Hence it is uniformly bounded. The same argument proves that $s_1 $ is uniformly bounded
on $K$ and concludes the proof of lemma \ref{lem14}.

\qed

\begin{lemma}\label{lem15}
Let $x\in \Tc (f)$ be such that its   orbit is relatively compact. Then we have~:
$\displaystyle{\lim_{n\rightarrow \infty} b_n(x)=+\infty}$.
\end{lemma}
Let us notice that it gives exactly the first part of proposition \ref{prop12}.

\noindent{\bf Proof of lemma \ref{lem15}~:} We will use a change of basis along the   orbit of $x$~:  let us denote by $s_-(f^kx)$ the slope of $G_-(f^kx)$ and by $s_+(f^kx)$ the slope of $G_+(f^kx)$. We will choose $G_-(x)$ as new``horizontal line'', i.e. if the ``old coordinates'' in $T_{y}\A$ are $(\Theta, R)$, the new coordinates are~:
$$P(y).\begin{pmatrix}
\Theta\\
R\\
\end{pmatrix}=\begin{pmatrix} 1&0\\
-s_-(y)&1\\
\end{pmatrix}
\begin{pmatrix}
\Theta\\
R\\
\end{pmatrix}=
\begin{pmatrix}
\Theta\\
-s_-(y)\Theta +R\\
\end{pmatrix}$$
In general, $P$ doesn't depend continuously on the considered point, but by lemma \ref{lem14}, $P$ and $P^{-1}$ are uniformly bounded along the orbit of $x$ (because $s_-$ is uniformly bounded). Moreover, $P$ is symplectic. Let us compute in the new coordinates $N_n(x_k)=P(x_{n+k})M_n(x_k)P(x_k)^{-1}$~:
$$N_n(x_k)=\begin{pmatrix}
b_n(x_k)(s_-(x_k)-s_{-n}(x_k))& b_n(x_k)\\
0&b_n(x_k)(s_n(x_{k+n})-s_+(x_{k+n}))\\
\end{pmatrix}.$$
We know that~: $\displaystyle{\lim_{n\rightarrow \infty}\uparrow s_{-n}(x_k)=s_-(x_k)}$. Hence~: $\displaystyle{\lim_{n\rightarrow +\infty}(s_-(x_k)-s_{-n}(x_k))=0^+}$. By lemma \ref{lem14}~: $\forall n\geq 1, s_n(x_{k+n})-s_+(x_{k+n})\leq 2A$. As $N_n$ is symplectic, we have~: $1=\det N_n(x_k)=b_n(x_k)^2(s_-(x_k)-s_{-n}(x_k))(s_n(x_{k+n})-s_+(x_{k+n})).$ We deduce~:
$$\forall n\in\N^*, 1\leq 2Ab_n(x_k)^2(s_-(x_k)-s_{-n}(x_k)) $$
and then~: $\displaystyle{\lim_{n\rightarrow \infty}b_n(x_k)=+\infty}$.

Let us now prove corollary \ref{cor13}. Let us assume that $v\in T_x\A\backslash G_-(x)$. We use the ``old coordinates''  (the usual ones) and write~: $v=(v_1,v_2)$. Because $v\notin G_-(x)$, we have~: $s_-(x)v_1-v_2\not=0$  and we compute~: $D\pi\circ Df^n(x)(v_1,v_2)= b_n(x)(v_2-s_{-n}(x)v_1)$ with $\displaystyle{\lim_{n\rightarrow +\infty} (v_2-s_{-n}(x)v_1)=v_2-s_-(x)v_1\not=0}$ and $\displaystyle{\lim_{n\rightarrow +\infty}b_n(x)=+\infty}$. We deduce that~:
$$\lim_{n\rightarrow +\infty} |D\pi\circ Df^n(x)v|=+\infty.$$\qed

\subsection{Some easy consequences concerning (non uniform) hyperbolicity}
All the results contained in this subsection are not new, see for example \cite{C-I}.
At first, an easy and well-known consequence of the dynamical criterion is the following~:

\begin{prop}\label{prop16} (Contreras-Iturriaga)
Let $M$ be an $f$ ordered and uniformly hyperbolic set where $f$ is an exact symplectic positive twist map. Then at every $x\in M$, $G_-(x)=E^s(x)$ and $G_+(x)=E^u(x)$ are transverse.
\end{prop}
The argument is only the characterization of the stable and unstable tangent spaces for an uniformly hyperbolic set and the dynamical criterion for $G_-$ and $G_+$.
\smallskip

Let us now consider an irrational  Mather measure $\mu$ for a positive twist map $f$. We have noticed that $\mu$ is  ergodic. Hence we can associate to $\mu$ two Lyapunov exponents, $-\lambda$ and $\lambda$ (because $f$ is area preserving). If $\lambda\not=0$, we say that the measure is (non uniformly) hyperbolic and  the Oseledet theorem asserts that at $\mu$ almost all points there exists a measurable splitting $T_x\A=E^s_x\oplus E^u_x$ in two transverse lines, invariant under $Df$ such that~:
\begin{enumerate}
\item[$\bullet$] $\displaystyle{\forall v\in E^s_x, \lim_{n\rightarrow +\infty} \| Df^{n}(x)v\|=0}$;
\item[$\bullet$]$\displaystyle{\forall v\in E^u_x, \lim_{n\rightarrow +\infty} \| Df^{-n}(x)v\|=0}$.
\end{enumerate}
Then we (classically) deduce from the dynamical criterion that~: $G_-(x)=E^s(x)$ and $G_+(x)=E^u(x)$ are $\mu$ almost everywhere transverse~:

\begin{prop}\label{prop17}(Contreras-Iturriaga)
Let $\mu$ be a Mather measure of an exact symplectic positive twist map. If the Lyapunov exponents of $\mu$ are non zero, then at $\mu$ almost all points, $G_-$ and $G_+$ are transverse.
\end{prop}

We have explained why, for  (non uniformly) hyperbolic Mather measures, the Green bundles are almost everywhere transverse. We will now interest ourselves in  the converse assertion~: if the Green bundles are (almost everywhere)transverse, is the dynamic (non uniformly) hyperbolic?\\
We  begin by the uniform case, and then consider the non uniform one. 

\subsection{What happens when the Green bundles are  everywhere transverse}  
It is known that, with some additional hypothesis, the transversality of the Green bundles implies hyperbolicity. For example in \cite{C-I}, the authors prove that if $K\subset \Tc(f)$ is an invariant compact subset such that on $K$, the Green bundles are transverse and such that $f_{|K}$ is non wandering, then $K$ is hyperbolic for $f$. As we know that the dynamic on Aubry-Mather sets is minimal and then non wandering, we can deduce a result for the Aubry-Mather sets. \\
In fact, we notice that the hypothesis ``$f_{|K}$ is non wandering'' is useless and that's why we give a new statement~:

\begin{thm}\label{thm18}  Let $f$ be an exact symplectic positive $C^1$ twist map and let $K\subset \Tc (f)$ be an invariant compact subset of $\Tc (f)$ such that, at every point of $K$, $G_-(x)$ and $G_+(x)$ are transverse. Then $K$ is uniformly hyperbolic and at every $x\in K$, we have~: $G_-(x)=E^s(x)$ and $G_+(x)=E^u(x)$.
\end{thm}

\begin{cor}\label{cor19}
Let $M$ be an irrational Aubry-Mather set for an exact symplectic positive $C^1$ twist map $f$ such that, at every point of $M$, $G_-(x)$ and $G_+(x)$ are transverse. Then $M$ is uniformly hyperbolic and at every $x\in M$, we have~: $G_-(x)=E^s(x)$ and $G_+(x)=E^u(x)$.
\end{cor}

\begin{cor}\label{cor20}
Let $M$ be an irrational Aubry-Mather set for an exact symplectic positive $C^1$ twist map $f$ which is not uniformly hyperbolic. There $\Gc (M)=\{ x\in M; G_-(x)=G_+(x)\}$ is a dense $G_\delta$-subset of $M$ and at every $x\in \Gc (M)$, $M$ is $C^1$ regular.
\end{cor}
This corollary is a consequence of theorem \ref{thm18} and proposition \ref{prop10}.
In order to prove theorem \ref{thm18}, let us give a definition~:

\begin{defin}
Let $(F_k)_{k\in\Z}$ be a continuous cocycle on a linear normed bundle $P~: E\rightarrow K$ above a compact metric space $K$. We say that the cocycle is {\em quasi-hyperbolic} if~: 
$$\forall v\in E, v\not=0\Rightarrow \sup_{k\in\Z}\| F_k v\|=+\infty.$$
\end{defin}
A consequence of the dynamical criterion (corollary \ref{cor13}) is~: if $K\subset \Tc (f)$ is a compact invariant subset of $\Tc (f)$ such that for every $x\in K$, $G_+(x)$ and $G_-(x)$ are transverse, then $(Df^k_{|K})_{k\in\Z}$ is a quasi-hyperbolic cocycle. Hence, we only have to prove the following statement to deduce the proof of theorem \ref{thm18}~:
\begin{thm}\label{thm21}
Let $(F_k)$ be a continuous, symplectic and quasi-hyperbolic cocycle on a linear and symplectic  (finite dimensional)bundle $P~: E\rightarrow K$ above a compact metric space $K$. Then $(F_k)_{k\in\Z}$ is hyperbolic.
\end{thm}

Let us give two lemmas which will be useful to prove this theorem.   The ideas of these lemmas are not new and the reader can find similar statements in the setting of the so-called ``quasi-Anosov diffeomorphisms'' for example in \cite{Man1}.

\begin{lemma}\label{lemm22}
Let $(F_k)_{k\in\Z}$ be a continuous and quasi-hyperbolic cocycle on a linear normed bundle $P~: E\rightarrow K$ above a compact metric space $K$. Let us define~:
\begin{enumerate}
\item[$\bullet$] $\displaystyle{E^s=\{ v\in E; \sup_{k\geq 0} \| F_kv\|<\infty\}}$;
\item[$\bullet$] $\displaystyle{E^u=\{ v\in E; \sup_{k\leq 0} \| F_kv\|<\infty\}}$.
\end{enumerate}
Then $(F_{n|E^s})_{n\geq 0}$ and $(F_{-n|E^u})_{n\geq 0}$ are uniformly contracting.
\end{lemma} 

\begin{lemma}\label{lemm23}
Let $(F_k)_{k\in\Z}$ be a continuous and quasi-hyperbolic cocycle on a linear normed bundle $P~: E\rightarrow K$ above a compact metric space $K$. If $(x_n)$ is a sequence of points of $K$ tending to $x$ and $(k_n)$ a sequence of integers tending to $+\infty$ such that $\displaystyle{\lim_{n\rightarrow \infty} P\circ F_{k_n}(x_n)=y\in K}$, then $\dim E^u(y)\geq codim E^s(x)$.
\end{lemma}

\noindent Let us explain how to deduce theorem \ref{thm21} from these lemmas~:
\smallskip

\noindent{\bf Proof of theorem \ref{thm21}~:} If the dimension of $E$ is $2d$, we only have to prove that~: $\forall x\in K, \dim E^u(x)=\dim E^s(x)=d$. Let us prove for example that $\dim E^u(x)=d$. \\
By lemma \ref{lemm22}, $(F_{n|E^s})_{n\geq 0}$ and $(F_{-n|E^u})_{n\geq 0}$ are uniformly contracting. As the cocycle is symplectic, we deduce that every $E^s(x)$ and $E^u(x)$ is isotropic for the symplectic form  and then $\dim E^s(x)\leq d$ and $\dim E^u(x)\leq d$.\\
 Let us now consider $x\in K$. As $K$ is compact, we can find a sequence $(k_n)_{n\in\N}$ of integers tending to $+\infty$ such that the sequence $(P\circ F_{k_n}(x))_{n\in\N}$ converges to a point $y\in K$. Then, by lemma \ref{lemm23}, we have~:  $\dim E^u(y)\geq codim E^s(x)$. But we know that $\dim E^u(y)\leq d$, hence~: $  2d-dim E^s(x)\leq \dim E^u(y)\leq d$ and~: $\dim E^s(x)=d$.
 \enddemo
 Let us now prove the two lemmas~:
 
 \noindent{\bf Proof of lemma \ref{lemm22}~:} We will only prove the result for $E^s$.
 
 Let us assume that we know that~:
 $$(*) \forall C>1, \exists N_C\geq 1, \forall v\in E^s, \forall n\geq N_C, \| F_nv\| \leq \frac{ \sup\{ \| F_kv\| ; k\geq 0\}}{C}.$$
 Then in this case~: $\sup \{ \| F_k v\| ; k\geq 0\} =\sup\{ \| F_k v\| ; k\in |[0, N_C]|\}$. We define~: $M=\sup\{ \| F_k (x)\| ; x\in K, k\in |[0, N_C]|\}$.   Then, if $j\in |[0, N_C-1]|$ and $n\in\N$~:
 $$\| F_{nN_c+j}v\| \leq \frac{1}{C} \sup\{ \|  F_{(n-1)N_C+j+k}v\|; k\geq 0\}\leq \frac{1}{C^2} \sup\{ \|  F_{(n-2)N_C+j+k}v\|; k\geq 0\}$$
 $$\dots \leq \frac{1}{C^n} \sup\{ \|  F_{j+k}v\|; k\geq 0\}\leq \frac{1}{C^n} \sup\{ \|  F_{k}v\|; k\geq 0\}\leq \frac{M}{C^n}\| v\|.$$
 This prove exponential contraction. \\
 
 Let us now prove $(*)$. If $(*)$ is not true, there exists $C>1$, a sequence $(k_n)$ in $\N$ tending to $+\infty$ and $v_n\in E^s$ with $\| v_n\|=1$ such that~:
 $$\forall n\in\N, \| F_{k_n}v_n\| \geq \frac{\sup\{ \| F_k v_n\|; k\geq 0\}}{C}.$$ 
 Then we define~: $w_n=\frac{F_{k_n}(v_n)}{\| F_{k_n}(v_n)\|}$. If we take a subsequence, we can assume that the sequence $(w_n)$ converges to a limit $w\in E$. Then we have~:
 $$\forall n\in\N, \forall k\in[-k_n, +\infty[, \| F_kw_n\|=\frac{\| F_{k+k_n}(v_n)\|}{\| F_{k_n}v_n\|}\leq \frac{\sup\{ \| F_jv_n\|; j\geq 0\} }{\| F_{k_n}v_n\|}\leq C.$$
 Hence~: $\forall k\in\Z, \| F_k w\|\leq C$; it is impossible because $\| w\|=1$ and the cocycle is quasi-hyperbolic.\enddemo

  \noindent{\bf Proof of lemma \ref{lemm23}~:} With the notation of this lemma, we choose a linear subspace $V\subset E_x$ such that $V$ is transverse to $E^s(x)$. What we want to prove is~: $\dim E^u(y)\geq \dim V$.\\
  We choose $V_n\in E_{x_n}$ such that $\displaystyle{\lim_{n\rightarrow \infty} V_n=V}$. If we use a subsequence, we have~: $\displaystyle{\lim_{n\rightarrow \infty} F_{k_n}(V_n)=V'\subset E_y}$. Then  we will prove~: $V'\subset E^u(y)$.\\
  
  Let us assume that we have proved that there exists $C>0$ such that~:
  $$(*) \forall n, \forall 0\leq k\leq k_n,  \| F_{-k|F_{k_n}(V_n)}\|\leq C.$$
  Then~: $\forall w\in V', \forall k\in \Z_-, \| F_kw\| \leq C\| w\|$ and $w\in E^u(y)$.
  
  Let us now assume that $(*)$ is not true~: we find $j_n\in \N$ and $i_n\in |[ 0, k_{j_n}||=|[0, K_n]|$ such that $ \| F_{-i_n|F_{K_n}(V_{i_n})}\|\geq n.$ If we extract a subsequence, we have $i_n\in |[0, k_n]|$ and  $\| F_{-i_n|F_{k_n}(V_n)}\|\geq n$. We choose $w_n\in F_{k_n}(V_n)$ such that~: $\| w_n\| =1$ and $\|F_{-i_n}(w_n)\| =\|  F_{-i_n|F_{k_n}(V_n)}\|$. We may even assume that~: $\| F_{-i_n}(w_n)\| = \sup\{ \| F_k(w_n)\| ; k\in |[-k_n, 0]|\} \geq n$.\\
  Then~: $\displaystyle{\lim_{n\rightarrow +\infty} i_n=+\infty}$. If $v_n=\frac{F_{-i_n}(w_n)}{\|F_{-i_n}(w_n)\|}$, we may extract a subsequence and assume that~: $\displaystyle{\lim_{n\rightarrow \infty} v_n=v}$.\\
   Then we have~: 
  $\forall k\in |[ 0, i_n]|, \| F_kv_n\| \leq \| v_n\|$ and then~: $\forall k\in\N, \| F_k v\|\leq \| v\|$ and $v\in E^s$.\\
Now, we have two cases~:
  \begin{enumerate}
  \item[$\bullet$] either $(k_n-i_n)$ doesn't tend to $+\infty$; we may extract a subsequence and assume that $\displaystyle{\lim_{n\rightarrow +\infty}(k_n-i_n)=N\geq 0}$; then~: $\displaystyle{  F_{-N}v=\lim_{n\rightarrow \infty}F_{i_n-k_n}(v_n)=\lim_{n\rightarrow \infty}\frac{F_{-k_n}(w_n)}{\|F_{-i_n}(w_n)\|}}$. We have~: $\frac{F_{-k_n}(w_n)}{\|F_{-i_n}(w_n)\|}\in V_n$ and then $F_{-N}v\in V$. Moreover, $F_{-N}v\in F_{-N}E^s=E^s$. 
  As $\| v\| =1$ and $V$ is transverse to $E^s_x$, we obtain a contradiction.
  \item[$\bullet$] or $\displaystyle{\lim_{n\rightarrow \infty} (k_n-i_n)=+\infty}$. Then we have~: 
  $$\forall k\in |[ -k_n+i_n, i_n]|, k-i_n\in |[-k_n, 0]|$$and~: $ \| F_kv_n\|=\frac{\| F_{k-i_n}(w_n)\|}{\| F_{-{i_n}}w_n\|} \leq 1= \| v_n\|.$ Hence~: $\forall k\in \Z_-, \| F_kv\| \leq \| v\|$ and $v\in E^s\cap E^u$. This contradicts $\| v\|=1$ and the fact that the cocycle is quasi-hyperbolic.
  \end{enumerate}
  \enddemo

\subsection{What happens for the Mather measures whose   Green bundles are almost everywhere  transverse}
Let us now consider a Mather measure of $f\in\Mc_\omega^+$. The map $d~: {\rm supp}\mu\rightarrow \{ 0, 1\}$ defined by $d(x)=\dim(G_-(x)\cap G_+(x))$ being measurable and constant along the orbits of $f$, we know that $d$ is $\mu$-almost everywhere constant. This constant is 0 or 1. In this subsection, we will study the case of a constant equal to zero and  prove~: 

\begin{thm}\label{thm24} Let $f\in \Mc_\omega^+$ be an exact symplectic positive twist map and let $\mu$ be an irrational Mather measure for $f$. We assume that at $\mu$-almost every point, $G_-$ is transverse to $G_+$. Then the Lyapunov exponents of $\mu$ are non zero.
\end{thm}

\begin{cor}\label{coro25}
Let $f\in \Mc_\omega^+$ be an exact symplectic positive twist map and let $\mu$ be an irrational Mather measure for $f$. We assume that the Lyapunov exponents of $\mu$ are  zero. Then $\mu$ almost everywhere, ${\rm supp}\mu$ is  $C^1$ regular.
\end{cor}
Indeed, in this case, $d=\dim (G_-\cap G_+)$ is $\mu$-almost equal to 1, i.e. $\mu$-almost everywhere we have~: $G_-=G_+$. Hence we deduce from proposition \ref{prop10} that $\mu$-almost everywhere, ${\rm supp}(\mu )$ is $C^1$-regular. We deduce~:
\begin{cor}\label{coro26}
Let $f\in \Mc_\omega^+$ be an exact symplectic positive twist map and let $\mu$ be an irrational Mather measure for $f$. We assume that $\mu$ almost everywhere, ${\rm supp}(\mu )$ is not $C^1$-regular. Then  the Lyapunov exponents of $\mu$ are  non zero. \end{cor}

\noindent{\bf Proof of theorem \ref{thm24}~:} We will use the same notations as in the proof of proposition \ref{prop12}. At $x\in {\rm supp}\mu$, we have~:
$$M_n(x)=\begin{pmatrix}
-b_n(x)s_{-n}(x)&b_n(x)\\
-b_n(x)^{-1}-b_n(x)s_{-n}(x)s_n(x_n)&s_n(x_n)b_n(x)\\
\end{pmatrix}
$$
 Instead of using a change of basis which sends $G_-$ on the horizontal, we will use such a change which sends $G_+$ on the horizontal~: 
 $$P(x)=\begin{pmatrix}
 1&0\\
 -s_+(x)&1\\
 \end{pmatrix}
 $$
 In the new coordinates, the new matrix of $Df^n(x)$ is $N_n(x)=P(x_n)M_n(x)P(x)^{-1}$ with~:
 $$ N_n(x)=\begin{pmatrix}
 b_n(x)(s_+(x)-s_{-n}(x))& b_n(x)\\
 0&b_n(x)(s_n(x_n)-s_+(x_n))\\
 \end{pmatrix}.
 $$
 We will use in the proof lemma \ref{lem14} and two other lemmas~:
 
 \begin{lemma}\label{lemm25} Let $\varepsilon >0$. There exists a subset $K_\varepsilon\subset{\rm supp}\mu$ such that $\mu(K_\varepsilon )>1-\varepsilon$ and such that on $K_\varepsilon$, $(s_{-n})$ and $(s_n)$ converge uniformly on $K_\varepsilon$ to their limits $s_-$ and $s_+$.    \end{lemma}
  This lemma is just a consequence of Egorov theorem (see for example \cite{Ko-Fo}). 
    
   \begin{lemma}\label{lemm26}  Let $\varepsilon >0$. There exists a subset $F_\varepsilon\subset{\rm supp}\mu$  such that $\mu(F_\varepsilon) >1-\varepsilon$ and and $\alpha>0$ such that ~: 
  $\forall x\in F_\varepsilon, s_+(x)-s_-(x)\geq \alpha$.
    \end{lemma}
\noindent{\bf Proof of lemma \ref{lemm26}~:} We have assume that at $\mu$-almost every point $x\in\A$, $G_-(x)$ and $G_+(x)$ are transverse, i.e. $s_+(x)-s_-(x)>0$.  Hence~:
$$\mu\left( \bigcup_{n\geq 1}\{ x; s_+(x)-s_-(x)\geq \frac{1}{n}\}\right)=1.$$
As the previous union is monotone, we deduce that there exists $n\geq 1$ such that~:
$$\mu\left( \{ x; s_+(x)-s_-(x)\geq \frac{1}{n}\}\right)\geq 1-\varepsilon.$$
\enddemo
We deduce from these two lemmas that there exists $J_\varepsilon$ and a constant $\alpha>0$ such that $\mu(J_\varepsilon)\geq 1-\varepsilon$, $(s_n)$ and $(s_{-n})$ converge uniformly on $J_\varepsilon $ and~: $\forall x\in J_\varepsilon, s_+(x)-s_-(x)\geq \alpha$.
\begin{lemma}\label{lemm27}
Let $A>0$ and $\varepsilon >0$. Then there exists $N=N(A, \varepsilon)$ such that~:
$$\forall n\geq N, \forall x\in J_\varepsilon, f^nx\in J_\varepsilon \Rightarrow b_n(x)\geq A.$$
\end{lemma}
\noindent{\bf Proof of lemma \ref{lemm27}~:} We use the matrix $N_n(x)$~:
$1=\det N_n(x)=b_n(x)^2(s_+(x)-s_{-n}(x))(s_n(x_n)-s_+(x_n))$ with $x_n=f^n(x)$. By lemma \ref{lem14}, there exists $B>0$ such that~: $\forall y\in  {\rm supp} \mu, \forall k\in \Z, -B\leq s_k(x)\leq B$. Then~: $\forall x\in {\rm supp}\mu, \forall n\in \N^*, 0<s_+(x)-s_{-n}(x)\leq 2B$. We deduce~:
$\forall x\in {\rm supp}\mu, \forall n\in\N^*, 1\leq 2B b_n(x)^2 (s_n(x_n)-s_+(x_n))$.\\

By definition of $J_\varepsilon$, we know that $s_{n}$ converge uniformly on $J_\varepsilon$ to $s_+$. Hence there exists $N\geq 1$ such that~: $\forall n\geq N, \forall y\in J_\varepsilon, 0<s_n(y)-s_+(y)<\frac{1}{2BA^2}$.

Let us now assume that $x, x_n=f^n(x)\in J_\varepsilon$. Then~: $1\leq 2B b_n(x)^2 (s_n(x_n)-s_+(x_n))\leq 2B b_n(x)^2\frac{1}{2BA^2}=\frac{b_n(x)^2}{A^2}$ and $b_n(x)\geq A$.\enddemo
  To a given $\varepsilon>0$ we have associated a set $J_\varepsilon\subset {\rm supp}(\mu )$ such that $\mu (J_\varepsilon)> 1-\varepsilon$, $(s_n)$ and $(s_{-n})$ converge uniformly on $J_\varepsilon$ to their limits and $\forall x\in J_\varepsilon, s_+(x)-s_-(x)\geq \alpha>0$. By lemma \ref{lemm27}, we find $N\geq 1$ such that~:
  $$\forall x\in J_{\varepsilon}, \forall n\geq N, f^n(x)\in J_{\varepsilon} \Rightarrow b_n(x)\geq \frac{2}{\alpha}.$$
  Let us notice that because $\mu$ is an irrational Mather measure, it is ergodic not only for $f$ but  for $f^N$ too (we don't say that in general an ergodic measure for $f$ is ergodic for $f^N$, but this is true for $f$ homeomorphism of the circle with a irrational rotation number). If we denote by $\sharp Y$ the cardinal of a set $Y$, we know by the ergodic theorem of Birkhoff (see e.g. \cite{Man2}) that for almost $x\in J_\varepsilon$~:
  $$\frac{1}{\ell }\sharp\{ 0\leq k\leq \ell -1; f^{kN}(x)\in J_\varepsilon\}\stackrel{\ell\rightarrow +\infty}{\longrightarrow} \mu(J_\varepsilon ) \geq 1-\varepsilon.$$
  We denote by $\lambda$, $-\lambda$ the Lyapunov exponents of $f$ (with $\lambda \geq 0$).
  
  Then $L_\varepsilon$ is the set of points of $J_\varepsilon$ such that~:
  \begin{enumerate}
  \item[$\bullet$] $\frac{1}{\ell }\sharp\{ 0\leq k\leq \ell -1; f^{kN}(x)\in J_\varepsilon\}\stackrel{\ell\rightarrow +\infty}{\longrightarrow} \mu(J_\varepsilon )$;
   \item[$\bullet$] $x$ is a regular point for $\mu$ i.e. at $x$ there exists a splitting of the tangent space $T_x\A$ corresponding to the Lyapunov exponents (see e.g. \cite{Man2}).
  \end{enumerate}
  Then $\mu(L_\varepsilon)=\mu(J_\varepsilon)\geq 1-\varepsilon$ and if $x\in L_\varepsilon$, we have~: $\displaystyle{\lim_{n\rightarrow +\infty}\frac{1}{n}\log\| Df^n(x)\|=\lambda}$.\\
   If $x\in L_\varepsilon$, we define~:
  $$n(\ell )=\sharp\{ 0\leq k\leq \ell -1; f^{kN}(x)\in J_\varepsilon\}$$
  and $0=k(1)<k(2)<\dots < k(n(\ell ))\leq \ell$ are such that $f^{k(i)N}x\in J_\varepsilon$.\\
  The chain rule of derivatives implies that for all $x\in L_\varepsilon$~:\\
  $Df^{k(n(\ell))N}(x)=$
  $$Df^{(k(n(\ell ))-k(n(\ell )-1))N}(f^{k(n(\ell )-1)N}x).Df^{(k(n(\ell )-1)-k(n(\ell )-2))N}(f^{k(n(\ell )-2)N}x)\dots Df^{k(1)N}(x).$$
  We write this equality for the matrices $N_k$ and specially for the terms $a_k$~:\\
  $b_{k(n(\ell))N}(x)(s_+(x)-s_{-k(n(\ell))N}(x))
 =$
 $$b_{(k(n(\ell ))-k(n(\ell )-1))N}(f^{k(n(\ell )-1)N}x)\Delta s_{(k(n(\ell ))-k(n(\ell )-1))N}(f^{k(n(\ell )-1)N}x)\dots b_{k(1)N}(x)\Delta s_{k(1)N}(x).$$
 where~: $\Delta s_n(x):=s_+(x)-s_{-n}(x)$.\\
 Let us notice that~:  $\| Df^{k(n(\ell ))N}(x)\|\geq b_{k(n(\ell))N}(x)(s_+(x)-s_{-k(n(\ell))N})=a_{k(n(\ell ))}$. \\
 Moreover, as for every $0\leq j\leq n(\ell)$, we have $f^{k(j)N}(x)\in J_\varepsilon$, we know  that for every $0\leq j\leq n(\ell)-1$~:  $b_{(k(j+1)-k(j))N}(f^{k(j)N}(x))\geq \frac{2}{\alpha}$ and that $\Delta s_{(k(j+1)-k(j))N}(f^{k(j)N}(x))> s_+((f^{k(j)N}(x)))-s_-(f^{k(j)N}(x))\geq \alpha$. We deduce~:
 $$\| Df^{k(n(\ell))N}(x)\|\geq b_{k(n(\ell))N}(x)(s_+(x)-s_{-k(n(\ell))N}(x))\geq (\frac{2}{\alpha}.\alpha)^{n(\ell)}=2^{n(\ell )}$$
 We deduce~:
 $$\frac{1}{k(n(\ell ))N}\log \| Df^{k(n(\ell))N}(x)\|\geq \frac{n(\ell)}{k(n(\ell ))N}\log 2.$$
 But we have~: $k(n(\ell ))\leq \ell$ then~: $\frac{1}{k(n(\ell ))N}\log \| Df^{k(n(\ell))N}(x)\|\geq \frac{n(\ell)}{\ell N}\log 2$.\\
  As $\displaystyle{\lim_{\ell\rightarrow +\infty}\frac{n(\ell )}{\ell}=\mu (J_\varepsilon)\geq 1-\varepsilon}$, we obtain~:
  $$\lambda=\lim_{\ell\rightarrow +\infty}\frac{1}{k(n(\ell ))N}\log \| Df^{k(n(\ell))N}(x)\|\geq \frac{1-\varepsilon}{N}\log 2>0;
  $$
  hence the Lyapunov exponents are non zero.\enddemo
  
  \section{The hyperbolic case~: proof of its irregularity}
  \subsection{Case of uniform hyperbolicity}
  \begin{prop}\label{prop30}
  Let $M$ be an uniformly hyperbolic  irrational Aubry-Mather set of an exact symplectic positive $C^1$ twist map $f$ of $\A$. Then at every $x\in M$, $M$ is not $C^1$ regular.
  \end{prop}
  
  \noindent{\bf Proof of proposition \ref{prop30}~:} At first, let us notice that such a $M$ cannot be a curve~: we proved in \cite{Arna1} that if the graph of a continuous map $\gamma~:\T\rightarrow \R$ is invariant by $f$, then   Lebesgue almost everywhere we have~: $G_-(t,\gamma (t))=G_+(t, \gamma(t))$, which contradicts proposition\ref{prop16} which asserts that $G_-=E^s$ and $G_+=E^u$.  Another argument is the fact, proved in \cite{MMP} , that $\pi (M)$ has zero Lebesgue measure.
  
Hence $M$ is a Cantor and the dynamic on $M$ is Lipschitz conjugate to the one of a Denjoy counter-example on its minimal invariant set. Then we consider two points $x\not= y$ of $M$ such that there exists an open  interval $I\subset \T$ whose ends are $\pi(x)$ and $\pi(y)$ and which doesn't meet $\pi(M)$~: $I\cap \pi(M)=\emptyset$.  We deduce from the dynamic of the Denjoy counter-examples (see \cite{He1}) that~:
\begin{enumerate}
\item[$\bullet$] the positive and negative orbits of $x$ and $y$ under $f$ are dense in $M$;
\item[$\bullet$] $\displaystyle{\lim_{n\rightarrow +\infty}d(f^nx,f^ny)=\lim_{n\rightarrow +\infty}d(f^{-n}x,f^{-n}y)=0}$.
\end{enumerate}
As $M$ is uniformly hyperbolic, we can define a local stable and unstable laminations on $M$ (see for example \cite{Shu}), $W^s_{\rm loc}$ and $W^u_{\rm loc}$. Then for $n$ big enough, $f^nx$ and $f^ny$ belongs to the same local stable leaf, and $f^{-n}x$ and $f^{-n}y$ belongs to the same local unstable leaf. Hence, because $\displaystyle{\lim_{n\rightarrow +\infty}d(f^nx,f^ny)=\lim_{n\rightarrow +\infty}d(f^{-n}x,f^{-n}y)=0}$, for $n$ big enough, the vector joining $f^nx$ to $f^ny$ (resp. $f^{-n}x$ to $f^{-n}y$) is  close  the stable bundle $E^s$ (resp. the unstable bundle $E^u$).

Let now $z\in M$ be any point. Then there exists two sequences $(i_n)$ and $(j_n)$ of integers which tends to $+\infty$ and are such that~:
$$\lim_{n\rightarrow +\infty} f^{i_n}x=\lim_{n\rightarrow +\infty} f^{i_n}y=\lim_{n\rightarrow +\infty} f^{-j_n}x=\lim_{n\rightarrow +\infty} f^{-j_n}y=z.$$
The direction of the ``vector'' joining $f^{i_n}x$ to $f^{i_n}y$ tends to $E^s(z)$ and the direction of the vector joining $f^{-j_n}x$ to $f^{-j_n}y$ tends to $E^u(z)$. Hence~: $E^u(z)\cup E^s(z)\subset P_M(z)$ and $M$ is not $C^1$-regular at $z$.\enddemo

\subsection{Case of non uniform hyperbolicity}
\begin{prop}\label{prop31}
Let $f\in\Mc_\omega$ be an exact  symplectic positive $C^1$ twist map and let $\mu$ be an irrational Mather measure of $f$ whose Lyapunov exponents are non zero. Then, at $\mu$ almost every point, ${\rm supp}\mu$ is not $C^1$ regular. 
\end{prop}
To prove this result, we will need some results concerning ergodic theory (see for example \cite{Po}); for us, every probability space $(X,\mu)$ will be such that $X$ is a metric compact space endowed with its Borel $\sigma$-algebra.

\begin{defin}
Let $(X,\mu)$ be a probability space, $T$ be a measure preserving transformation of $(X,\mu)$    and $(f_n)\in L^1(X,\mu )$ be a sequence of $\mu$-integrable functions from $X$ to $\R$. Then $(f_n)$ is {\em $T$-subadditive} if for $\mu$ almost every $x\in X$ and all $n,m\in\N$, we have~: $f_{n+m}(x)\leq f_n(x)+f_m(T^nx)$.
\end{defin}
A useful result in ergodic theory is the following~:
\begin{prop}\label{prop32}  (Subadditive ergodic theorem, Klingman) Let $(X,\mu)$ be a probability space, let $T$ be a measure preserving transformation of $(X,\mu)$  such that $\mu$ is ergodic for $T$    and let $f=(f_n)\in L^1(X,\mu )$ be a $T$-subadditive sequence. Then there exists a constant $\Lambda (f)\geq -\infty$ such that for $\mu$-almost every $x\in X$, we have~:
$$\lim_{n\rightarrow +\infty} \frac{1}{n}f_n(x)=\Lambda (f).$$
Moreover, the constant $\Lambda (f)$ satisfies~:
$$\Lambda (f)=\lim_{n\rightarrow \infty} \frac{1}{n}\int f_n d\mu=\inf_n \frac{1}{n}\int f_n d\mu.$$

\end{prop}

We will use the following refinement of this proposition, which concerns only the uniquely ergodic measures. A proof of it in the case of continuous functions is given in \cite{Fur}; the proof for upper  semi-continuous functions is exactly the same. 
\begin{prop}\label{Pfur}
Let $(X,\mu)$ be a probability space, $T$ be a measure preserving transformation of $(X,\mu)$  such that $\mu$ is uniquely ergodic for $T$    and $(f_n)\in L^1(X,\mu )$ be a $T$-subadditive sequence of upper semi-continuous functions.  Let $\Lambda (f)$ be the constant associated to $f$ via the subadditive ergodic theorem. We assume that $\Lambda (f)\in\R$. Then~:
$$\forall \varepsilon >0, \exists N\geq 0, \forall n\geq N, \forall x\in X, \frac{1}{n}f_n(x)\leq \Lambda (f)+\varepsilon.$$
\end{prop}

\noindent{\bf Proof of proposition \ref{prop31}~:} At first, let us notice that the set $R$ of points where ${\rm supp}\mu$ is $C^1$ regular is a $G_\delta$ subset of ${\rm supp} \mu$ and then is measurable.  Let us assume that $\mu (R)=a>0$. If ${\rm supp}\mu$ is the graph of $\gamma$ above $\pi ({\rm Supp}\mu)$ then $\gamma$ is differentiable at every $\theta\in \pi(R)$ and even $C^1$ at such a $\theta$. Moreover, $R$ is invariant by $f$.\\\
We know that there exists an orientation preserving bi-Lipschitz homeomorphism $h~: \T\rightarrow \T$ such that for all $(\theta ,r)\in {\rm supp} \mu$, we have~: $\pi\circ f(\theta, r)=h(\theta)$. We denote by $m$ the unique $h$-invariant probability measure on $\T$ (this measure is supported in $\pi ({\rm supp} \mu)$). \\
We may choose $h$ in a more precise way~:  If $I=]a, b[$ is an open interval which is a connected component of $\T\backslash \pi ({\rm supp} \mu)$, we may choose $h$ affine on $I$. Let $D$ be the (countable) set of the points of $\pi ({\rm supp} \mu)$ which are ends of such intervals. Let us prove that every $h^k$ is differentiable on $\pi (R)\backslash D$~:\\
Let us consider $\theta\in\pi (R)\backslash D$ and $(\alpha_n)< (\beta_n)$ two sequences of elements of $\T$ converging to $\theta$. Let $I_n=[\alpha_n^1, \alpha_n^2]$ (resp. $J_n=[\beta_n^1, \beta_n^2]$) be~:
\begin{enumerate}
\item[$\bullet$]  either the longest closed interval of $(\T\backslash \pi ({\rm supp}\mu ))\cup D$ containing $\alpha_n$ (resp. $\beta_n$) if $\alpha_n\notin \pi({\rm supp}\mu)\backslash D$ (resp. $\beta_n\notin \pi({\rm supp}\mu)\backslash D$);
\item[$\bullet$] or $\{ \alpha_n\}$ (resp $\{\beta_n\}$) if $\alpha_n\in \pi({\rm supp}\mu)\backslash D$ (resp. $\beta_n\in \pi({\rm supp}\mu)\backslash D$).
\end{enumerate}
As $\theta\notin D$, we have~: 
$$\lim_{n\rightarrow \infty}\alpha_n^1=\lim_{n\rightarrow \infty}\alpha_n^2=\lim_{n\rightarrow \infty}\beta_n^1=\lim_{n\rightarrow \infty}\beta_n^2=\theta.$$
Moreover (we denote by CH the convex hull)~: 
$$\frac{h^k(\alpha_n)-h^k(\beta_n)}{\alpha_n-\beta_n}\in{\rm CH}\left\{ \frac{h^k(\alpha_n)-h^k(\alpha^2_n)}{\alpha_n-\alpha^2_n},\frac{h^k(\alpha^2_n)-h^k(\beta^1_n)}{\alpha^2_n-\beta^1_n},\frac{h^k(\beta^1_n)-h^k(\beta_n)}{\beta^1_n-\beta_n}\right\} .
$$
(when the written slope is not defined, we don't write it)\\
As $h^k$ is affine on $I_n$ and $J_n$,  this last set is equal to~: 
$${\rm CH}\left\{ \frac{h^k(\alpha^1_n)-h^k(\alpha^2_n)}{\alpha^1_n-\alpha^2_n},\frac{h^k(\alpha^2_n)-h^k(\beta^1_n)}{\alpha^2_n-\beta^1_n},\frac{h^k(\beta^1_n)-h^k(\beta^2_n)}{\beta^1_n-\beta^2_n}\right\}
$$
As $\alpha_n^1, \alpha^2_n, \beta^1_n, \beta^2_n\in \pi({\rm supp}\mu )$ tend to $\theta\in\pi(R)$ when $n$ goes to $+\infty$, we have (when the slope is defined i.e. $\alpha_n^1\not=\alpha_n^2$)~:\\
$$\lim_{n\rightarrow \infty} \frac{h^k(\alpha^1_n)-h^k(\alpha^2_n)}{\alpha^1_n-\alpha^2_n}=\lim_{n\rightarrow \infty}  \frac{\pi\circ f^k(\alpha^1_n,\gamma (\alpha^1_n))-\pi\circ f^k(\alpha^2_n,\gamma (\alpha^2_n))}{\alpha^1_n-\alpha^2_n}$$
$$=D\pi\circ Df^k(\theta, \gamma(\theta))(1, \gamma'(\theta ))$$
and similarly (if defined)~:
$$\lim_{n\rightarrow \infty} \frac{h^k(\alpha^2_n)-h^k(\beta^1_n)}{\alpha^2_n-\beta^1_n}=\lim_{n\rightarrow \infty} \frac{h^k(\beta^1_n)-h^k(\beta^2_n)}{\beta^1_n-\beta^2_n}=D\pi\circ Df^k(\theta, \gamma(\theta))(1, \gamma'(\theta ))$$
Hence~: $\displaystyle{ \lim_{n\rightarrow \infty} \frac{h^k(\alpha_n)-h^k(\beta_n)}{\alpha_n-\beta_n}
=D\pi\circ Df^k(\theta, \gamma(\theta))(1, \gamma'(\theta ))
}$.\\
Finally,  every $h^n$ is differentiable on $\pi (R)\backslash D$ and~: 
$$\forall \theta\in\pi(R)\backslash D, \forall n\in\N, \lim_{\alpha, \beta\rightarrow \theta}\frac{h^n(\alpha)-h^n(\beta)}{\alpha-\beta}=(h^n)'(\theta)=D\pi\circ Df^n(\theta, \gamma(\theta))(1, \gamma'(\theta )).$$
 
We define for every $\theta\in \T$~: $\displaystyle{h'_n(\theta)=\liminf_{y\not=z\rightarrow
\theta}\frac{h^n(z)-h^n(y)}{z-y}>0}$; then every $h_n'$
is lower semi-continuous and then measurable. As $h$
is bi-Lipschitz, there exists $K_n>1$ such that for
every $x\in\T$, $\frac{1}{K_n}\leq h'_n(x)\leq K_n$.
Hence every $g_n=-\log h_n'$ is bounded and measurable
and thus belongs to  $L^1(m)$ and  the sequence
$g=(g_n)_{n\geq 1}$ is a $h$-subadditive sequence.
Moreover, every $g_n$ is upper semicontinuous. As $m$
is uniquely ergodic for $h$, we may apply proposition
\ref{Pfur}~: 
$$\forall \varepsilon >0, \exists N\geq 0,\forall \theta\in\T,  \forall n\geq N, \frac{1}{n}g_n(\theta )\leq \Lambda (g)+\varepsilon.$$
Let $\lambda$ be the Lebesgue measure on $\T$. As $(-\log )$ is convex, we have by Jensen inequality~:
$$-\log\left( \int h'_n d\lambda \right)\leq -\int  \log  h'_n d\lambda= \int g_n d\lambda.$$
Moreover, $h$ being Lipschitz is $\lambda$-almost everywhere differentiable and~:
$\int h_n'd\lambda\leq \int(h^n)'d\lambda =\tilde h^n (1)-\tilde h^n(0)=1$. Hence~:
$$0=-\log 1\leq -\log\left( \int h'_n d\lambda \right)\leq   \int g_n d\lambda$$
i.e~: $\int g_nd\lambda \geq 0$. Let us now choose $\varepsilon >0$. We know that there exists $N\geq 1$ such that~:  $\forall n\geq N, \forall x\in\T, \frac{1}{n}g_n(x)\leq \Lambda (g)+\varepsilon$ and thus~: $\forall n\geq N, 0\leq \frac{1}{n}\int g_nd\lambda\leq \Lambda (g)+\varepsilon$. We deduce that~: $\Lambda (g)\geq 0$.\\
By proposition \ref{prop32},  we know that for $m$-almost $\theta\in \T$, we have~: 
$\displaystyle{\lim_{n\rightarrow +\infty} \frac{1}{n}g_n(\theta )=\Lambda (g)}.$ Hence for $m$-almost $\theta \in \pi(R)\backslash D$, we have~: $\displaystyle{\lim_{n\rightarrow +\infty} \frac{1}{n}g_n(\theta )\geq 0}$; we denote by $A=\pi(R')$ the set of such $\theta$. We have noticed that for such a $\theta$, if $(\theta, r)\in {\rm supp} \mu$~:
\begin{enumerate}
\item[$\bullet$] every $h^n$ is differentiable at $\theta$ and even~: $\displaystyle{(h^n)'(\theta)=\lim_{y\not=z\rightarrow \theta}\frac{h^n(z)-h^n(y)}{y-z} }=h_n(\theta) $ and then $g_n(\theta)=-\log ((h^n)'(\theta ))$;
\item[$\bullet$] we have seen too that~: $
(h^n)'(\theta)=D\pi\circ Df^n(\theta, r)(1, \gamma'(\theta ))$.
\end{enumerate}
Let us now denote by $\nu> 0, -\nu$ the Lyapunov exponents of $\mu$ for $f$. Then there exists a subset $S$ of $R'$ such that $\mu (S)=\mu(R') =a>0$ and such that  at every $(\theta, r)\in S$ we can define the Oseledet's splitting $E^s\oplus E^u$~ :
$$\forall v\in E^u(\theta, r), \lim_{n\rightarrow\pm \infty}\frac{1}{n}\log\| Df^nv\|=\nu;\quad
\forall v\in E^s(\theta, r), \lim_{n\rightarrow\pm \infty}\frac{1}{n}\log\| Df^nv\|=-\nu.$$
Then for $(\theta, r)\in S$, we have (we recall that $\gamma'$ is bounded because ${\rm supp}\mu$ is Lipschitz)~:
$$\forall n\in\N^*, \frac{1}{n}\log\| Df^n(\theta ,r)(1, \gamma'(\theta))\|=\frac{1}{n}\log\| ( (h^n)'(\theta),\gamma'(h^n(\theta))(h^n)'(\theta ))\|= $$
 $$ \frac{1}{n}\log |(h^n)'(\theta)| 
  +\frac{1}{n}\log \| (1, \gamma' (h^n(\theta )))\|=- \frac{1}{n}g_n(\theta)+\frac{1}{n}\log \| (1, \gamma' (h^n(\theta )))\|\stackrel{n\rightarrow \infty}{\longrightarrow }- \Lambda (g)\leq 0.$$
We deduce that $(1, \gamma'(\theta))\in E^s(\theta, r)$. A similar argument for $n$ going to $-\infty$ (replacing $f$ by $f^{-1}$ and $h$ by $h^{-1}$) proves that $(1, \gamma'(\theta)) \in E^u(\theta, r)$. As $E^u(\theta,r)\cap E^s(\theta, r)=\{ 0\}$, we obtain a contradiction.\enddemo

\section{Proof of the results contained in the  introduction}

\noindent{\bf Proof of theorem \ref{th1}~:} we assume that $\mu$ is an irrational Mather measure of $f\in\Mc_\omega$; considering $f^{-1}$ instead of $f$, we may assume that $f\in\Mc_\omega^+$.\\
1) Let us assume that for $\mu$-almost $x$, ${\rm supp}\mu$ is $C^1$-regular at $x$. Then by proposition \ref{prop31}, the Lyapunov exponents of $f$ are zero.\\
2) Let us assume that the Lyapunov exponents of $\mu$ are zero. Then we deduce from corollary \ref{coro25} that ${\rm supp}\mu$ is $C^1$-regular $\mu$-almost everywhere.\enddemo

\noindent{\bf Proof of proposition \ref{prop2}~:} we assume that $\mu$ is an irrational Mather measure of $f\in\Mc_\omega$; considering $f^{-1}$ instead of $f$, we may assume that $f\in\Mc_\omega^+$.\\
1)  Let us assume that for $\mu$-almost $x$, ${\rm supp}\mu$ is not $C^1$-regular at $x$. Then by theorem \ref{th1}, the Lyapunov exponents of $\mu$ are non zero.\\
2) Let us assume that the Lyapunov exponents of $\mu$ are non zero. Then by proposition \ref{prop31}, ${\rm supp}\mu$ is $C^1$-irregular at $\mu$-almost every point.\enddemo

\noindent{\bf Proof of proposition \ref{prop4}~:} we assume that $M$ is an irrational Aubry-Mather set of $f\in\Mc_\omega$; considering $f^{-1}$ instead of $f$, we may assume that $f\in\Mc_\omega^+$.\\
1) we assume that $M$ is nowhere $C^1$-regular. By proposition \ref{prop10}, at every $x\in M$, $G_+(x)$ and $G_-(x)$ are transverse. Hence by corollary \ref{cor19}, $M$ is uniformly hyperbolic.\\
2) we assume that $M$ is uniformly hyperbolic. Then by proposition \ref{prop30}, $M$ is nowhere $C^1$-regular.\enddemo

\noindent{\bf Proof of proposition \ref{prop5}~:} Let $f\in\Mc_\omega^+$ be an exact  symplectic twist map and let $\mu$ be an irrational  Mather measure of $f$ which is non uniformly hyperbolic, i.e. the Lyapunov exponents are non zero but the corresponding Aubry-Mather set $M={\rm supp}\mu$ is not uniformly hyperbolic. The set $\Gc$ of the points $x$ of $M$ where $G_-(x)=G_+(x)$ is a $G_\delta$ of $M$ which is invariant by $f$. As $f_{|M}$ is minimal, either $\Gc$ is empty or it is a dense $G_\delta$ of $M$. Moreover, by proposition \ref{prop10}, at every point of $\Gc$, $M$ is $C^1$-regular. \\
Hence we only have to prove that $\Gc\not=\emptyset$. By theorem \ref{thm18}, as $M$ is not uniformly hyperbolic, $\Gc\not=\emptyset$.\enddemo

\noindent{\bf Proof of theorem \ref{th6}~:}
Let $f\in\Mc_\omega^+$ be an exact symplectic twist map and let $\Cc$ be a $C^1$ invariant curve which is a graph such that $f_{|\Cc}$ is $C^1$ conjugate to a rotation. Then we know (see \cite{Arna1}, it is an easy consequence of the dynamical criterion) that at every $x\in\Cc$, $G_-(x)=G_+(x)$. \\
Then, by proposition \ref{prop10}, the map $(x\in\Tc(f)\rightarrow G_-(x))$ and $(x\in\Tc(f)\rightarrow G_+(x))$ are continuous at every point of $\Cc$.\\
Let $W$ be a neigbourhood  of $T^1\Cc$, the unitary tangent bundle to $\Cc$  in  $T^1\A$, the unitary tangent bundle to $\A$. We may assume that $W$ is ``symmetrically fibered convex'' (i.e. if $u, v\in W\cap T_x\A$, if $\R u\preceq \R w\preceq\R v$, then $w\in W$). Then there exists a neighbourhood $V$ of $\Cc$ in $\A$ such that for every $x\in\Tc(f)\cap V$, $G^1_-(x)$ and 
$G^1_+(x)$ are in $W$ where $G_-^1$ and $G_+^1$ refer to the unitary Green bundles. Hence for every Aubry-Mather set $M$ for $f$ contained in $V$~:
$\forall x\in M, G_-^1(x), G_+^1(x)\in W$.
\\
Moreover, we know by proposition \ref{prop10} that~: $G_-(x)\preceq P_M(x)\preceq G_+(x)$. We deduce that 
for every Aubry-Mather set $M$ for $f$ contained in $V$~:
$$\forall x\in M, P^1_M(x)\subset W.$$
 \enddemo

\end{document}